\input amstex 
\documentstyle{amsppt}
\input bull-ppt
\keyedby{bull437e/pah}

\let\sl\it 

\topmatter
\cvol{30}
\cvolyear{1994}
\cmonth{January}
\cyear{1994}
\cvolno{1}
\cpgs{70-75}
\ratitle
\title  Pseudo-periodic Homeomorphisms and 
\\ Degeneration 
 of Riemann Surfaces
\endtitle 
\shorttitle{Pseudoperiodic homeomorphisms and degeneration 
of Riemann surfaces}
\author Yukio Matsumoto And Jos\'e Mar\'ia 
Montesinos-amilibia \endauthor 
\shortauthor{Yukio Matsumoto and J. M. Montesinos-amilibia}
\address Department of Mathematical Sciences, University 
of Tokyo, 
Hongo, Tokyo 113, Japan\endaddress 
\address Departamento de Geometr\'\i a y Topolog\'\i a, 
Facultad de Ciencias
Matem\'aticas, Universidad Complutense de Madrid, 28040 
Madrid, Spain\endaddress   
\date July 13, 1992\enddate
\thanks The first author was supported by DGICYT and 
Grant-in-Aid for 
Scientific Research
(C)03640020\endthanks
\thanks The second author was supported by DGICYT 
PB-89-0105\endthanks
\abstract We will announce two theorems. The first theorem 
will classify
all topological types of degenerate fibers appearing in 
one-parameter
families of Riemann surfaces, in terms of 
``pseudoperiodic'' surface
homeomorphisms. The second theorem will give a complete 
set of conjugacy
invariants for the mapping classes of such homeomorphisms. 
This latter 
result implies that Nielsen's set of invariants [{\it 
Surface
transformation classes of algebraically finite type}, 
Collected
Papers 2, Birkh\"auser (1986)] is not 
complete.\endabstract
\endtopmatter

\document

Let $\{ F_{\xi} \}$ be a family of Riemann surfaces 
parametrized by complex
numbers ${\xi}$. As ${\xi}$ approaches a special value, 
say, $0$, $F_{\xi}$
changes its ``shape'' and finally gets singularities 
becoming 
a singular surface $F_0$. This degeneration phenomenon has 
long been studied. Here we study
it from the topological point of view. We will show that 
the topological types
of the degenerate fibers can be completely classified in 
terms of certain
surface mapping classes introduced by Nielsen [Ni2] some 
fifty years 
ago. We will also give a complete set of conjugacy 
invariants for such 
mapping classes. \par
 
Throughout this paper all manifolds will be oriented, and 
all homeomorphisms 
between them will be orientation-preserving.  ${\Sigma}_g$ 
will denote 
a closed surface of genus $g$. Details will appear in [MM2].

\heading 1. Pseudoperiodic homeomorphisms \endheading
A homeomorphism $f : {\Sigma}_g \to {\Sigma}_g $ and its 
mapping class $[f\,]$
are called in this paper {\sl {pseudoperiodic}} if $[f\,]$ 
is either of finite 
order or reducible and in the latter case all component 
mapping classes
are of finite order. (Cf. [Th, G].) It was Nielsen [Ni2] 
who first 
studied these mapping classes under the name of 
{\sl{surface transformation
classes of algebraically finite type}}. Let us recall some 
conjugacy invariants
introduced by Nielsen. (See [Ni2, G].)\par

Suppose $f$ is reduced by a system of simple closed curves 
$C=C_1 \cup C_2 \cup
\cdots \cup C_r $. $C$ is called {\sl{admissible}} if each 
connected 
component of ${\Sigma}_g - C$ has negative Euler 
characteristic. If $g \ge 2$,
such a system always exists. With each curve $C_j$ of an 
admissible system $C$
is associated a rational number $s(C_j)$ called the {\sl{ 
screw number}}.
This measures the amount of Dehn twist performed by 
$f^{\alpha}$ about $C_j$,
where $\alpha = \alpha (C_j)$ is the smallest positive 
integer such that 
$f^{\alpha}(\overrightarrow{C_j}) = \overrightarrow{C_j}$. 
 An admissible 
system $C$ is {\sl{precise}} if $s(C_j) \neq 0$ for each 
$C_j$. A precise system always exists and is unique up to 
isotopy. 

A curve $C_j$ is {\sl{amphidrome}} if $\alpha$ is even and 
$f^{\alpha/2}
(\overrightarrow{C_j}) = -\overrightarrow{C_j}$. \par

We say that a pseudoperiodic homeomorphism $f$ is {\sl{of 
negative twist}}
if either $[f\,]$ is of finite order, or, when $[f\,]$ is 
reducible,
$s(C_j) < 0$ for each curve $C_j$ in a precise system $C$.

\heading 2. Degenerating family \endheading
By a {\sl{degenerating family \RM(of Riemann surfaces\RM) 
of genus $g$}} we 
mean a triple $(M,D, \varphi)$ consisting of a noncompact 
complex 
surface $M$; an open unit disk $D=\{ \xi \in \Bbb C \mid 
\vert \xi \vert < 1
\} $; and a surjective, proper, and holomorphic map $ 
\varphi : M \to D $. 
All fibers $F_{\xi} = {\varphi}^{-1}(\xi)$ are assumed to 
be connected,
and outside the origin ${\varphi} \vert  
{\varphi}^{-1}(D^{*}) :
{\varphi}^{-1} (D^{*}) \to D^{*}$ is
assumed to be a smooth fiber bundle with fiber 
${\Sigma}_g$ , where 
$D^{*} = D -\{ 0 \}$. The family is {\sl{minimal}} if it 
is free from 
$(-1)$-curves. Two families, $(M_i,D_i, {\varphi}_i), 
i=1,2$, are 
{\sl{topologically equivalent}} $(\overset\text{TOP}\to 
\sim)$ if there exist 
homeomorphisms $H:M_1 \to M_2$ and $h:
 D_1 \to D_2$ satisfying $h(0) = 0$
and $h{\varphi}_1 = {\varphi}_2 H$. We are interested in 
the following set:

$$ {\Cal S}_g = \{ \text{minimal degenerating families of 
genus } g \} / 
{\overset\text{TOP} \to \sim}.$$

Given a degenerating family of genus $g$, the {\sl 
{monodromy homeomorphism}}
$f:{\Sigma}_g \to {\Sigma}_g $ around the central fiber 
$F_0$ is determined
as usual (up to isotopy and conjugation). By the results 
of Imayoshi [I],
Shiga and Tanigawa [ST], and Earle and Sipe [ES], $f$ is a 
pseudoperiodic homeomorphism of
 negative twist. (There is an alternative topological 
proof, [MM2].)
 Let ${\Cal M}_g$ be the mapping class group of 
${\Sigma}_g$ and
 $\widehat{\Cal M}_g$ 
the set of conjugacy classes
in ${\Cal M}_g$. Let ${\Cal P}_g^{-}$ denote the subset of 
$\widehat{\Cal M}_g$
represented by pseudoperiodic mapping classes of negative 
twist. Then we have 
a well-defined map 

$$
\text{monodromy }\rho : {\Cal S}_g \to {\Cal P}_g^{-}.$$

\proclaim{Theorem 1} For $g \ge 2$, $\rho : {\Cal S}_g \to 
{\Cal P}_g^{-}$ 
is bijective. \endproclaim 

The corresponding map for $g=1$ is surjective, and the 
``kernel'' consists
of multiple fibers. (Cf. [K].) 
Using Theorem 1 and the construction in \S 3, one can 
topologically recover Namikawa and Ueno's classification of 
singular fibers of genus 2 [NU]. 
Since ${\Cal M}_g \cong \roman{Aut}({\pi}_1 {\Sigma}_g ) / 
\roman{Inn} ({\pi}_1 {\Sigma}_g )$,
we have 

\proclaim{Corollary 1.1} If $g \ge 2$, the action of the 
monodromy on 
${\pi}_1 {\Sigma}_g$ determines the topological 
equivalence class of 
$(M, D,\varphi)$. In particular, if the action is trivial, 
 $F_0$ is 
nonsingular. \endproclaim

Note that the action on $H_1 ( {\Sigma}_g ; \Bbb Z)$ is 
not sufficient [NU].
For an  explicit algebraic calculation of nonabelian 
monodromy see [O].

\proclaim{Corollary 1.2} Given a pseudoperiodic 
homeomorphism of 
negative twist $f :{\Sigma}_g \to {\Sigma}_g$, there 
exists a 
degenerating family $(M,D, \varphi)$ whose monodromy 
homeomorphism coincides
with $f$ up to isotopy and conjugation. \endproclaim

A closely related existence theorem has been independently 
announced by 
Earle and Sipe [ES, \S 7]. See [MM1] for a short  abstract 
of our result,
where we adopted a sign convention opposite to the one here.

\heading 3. Generalized quotient \endheading
The idea in proving Theorem 1 is to construct the inverse 
map of 
$\rho :
 {\Cal S}_g \to {\Cal P}_g^{-}$ .
A {\sl{Riemann surface with nodes}} $S$ was introduced by 
Bers [B].
We will call the underlying topological space of $S$ a 
{\sl{chorizo
space}} (chorizo = Spanish sausage), which we allow to 
have boundaries 
and not to be connected.
A chorizo space below will be {\sl{numerical}} in the 
sense that to 
each irreducible component is attached a positive integer 
called the
{\sl{multiplicity}}. \par

For a pseudoperiodic homeomorphism of negative twist 
$f:{\Sigma}_g \to {\Sigma}_g $ we can construct a 
numerical chorizo space
called the {\sl{generalized quotient}} $S_f$ of $f$ as 
follows: \par

Decompose ${\Sigma}_g$ as ${\Sigma}_g = A \cup B $, where 
$A$ is 
the union of annular neighborhoods of the curves in 
the precise system $C$ such that $f(A) =A$.
We assume that $f \vert B : B \to B $ is periodic. The 
quotient space 
$B/(f\vert B) $ is an orbifold. Let $p$ be a cone point, 
$(m, \lambda,
\sigma)$ the {\sl{valency}} of $p$ [Ni1]; that is, if $x 
\in B$ is a 
point over $p$, $m$ is the smallest positive integer such 
that 
$f^m(x)=x$, $f^m$ is the rotation around $x$ through the 
angle
$2\pi\delta/\lambda \; ( 0<\delta<\lambda, \roman{gcd}
(\lambda,\delta)=1)$ and
$\sigma$ is the integer determined by $\delta\sigma \equiv 
1 \; (\roman{mod} 
 \; \lambda)$, $0<\sigma < \lambda$.
\par

By the Euclidean algorithm we obtain a sequence of integers
$n_0>n_1>\cdots>n_l=1$ such that $n_0=\lambda, n_1 = 
\lambda-\sigma,
n_{i-1}+n_{i+1}\equiv 0 \; (\roman{mod} 
\; n_i), i=1, \dotsc, l-1$. Set $m_i = m{n_i}
 \; (i=0,1, \dotsc, l)$. 
\par
Let $\roman{Ch}(B)$ be the chorizo space constructed from 
$B/(f\vert B)$
by replacing a neighborhood of each cone point with the 
numerical
chorizo space shown in Figure
1, which consists of a disk and $l$ spheres.

Let $A_j (\subset A)$ be an annular neighborhood of $C_j$. 
The boundary 
curves $S_1$ and $S_2$ of $A_j$ have their 
{\sl{valencies}} $(m^{(1)},
{\lambda}^{(1)},{\sigma}^{(1)})$ and 
$(m^{(2)},{\lambda}^{(2)},{\sigma}^{(2)}
)$, when regarded as boundary curves of the periodic part 
$B$ [Ni1].
\par
Suppose $C_j$ is {\sl{not}} amphidrome.
Then $m^{(1)}=m^{(2)}={\alpha}(C_j)$.
Let $m$ be this common value.

\proclaim{Lemma} There exists uniquely a sequence of 
positive integers
$n_0, n_1, \dotsc , n_l \; (l\ge 1)$ satisfying the 
following conditions\RM:
\roster
\item"{\rm (i)}"
$n_0= {\lambda}^{(1)}, n_l={\lambda}^{(2)};$ 
\item"{\rm (ii)}"
$n_1 \equiv {\sigma}
^{(1)} (\roman{mod} \; {\lambda}^{(1)}), n_{l-1} \equiv 
{\sigma}^{(2)}
\; (\roman{mod} \; {\lambda}^
{(2)});$ 
\item"{\rm (iii)}" $n_{i-1} +n_{i+1} \equiv  0 \; 
(\roman{mod} \; n_i), i=1,2, \dotsc,
l-1; $
\item"{\rm (iv)}" $(n_{i-1}+n_{i+1})/n_i \ge 2, 
i=1,2,\dotsc, l-1;$ and
\item"{\rm (v)}" $\sum_{i=0}^{l-1} 1/n_i n_{i+1} = \vert 
s(C_j) \vert $. 
\endroster\endproclaim

Let $\roman{Ch}(A_j)$ be the chorizo space shown in Figure
2, which consists of 
two disks and $l-1$ spheres, $m_i$ being defined to be 
$mn_i$.

\fighere{5pc}
\caption{{\smc Figure 1}}

\topspace{5pc}
\caption{{\smc Figure 2}}

We consider the spaces $\roman{Ch}(f^iA_j)\ 
(i=0,1,\dotsc,m-1)$ identical : 
$\roman{Ch}(A_j)=
\roman{Ch}(fA_j)= \cdots =\roman{Ch}(f^{m-1}A_j)$. \par

Finally, suppose $A_j$ is amphidrome. Then $S_1$ and $S_2$ 
have
the same valency $(2m, \lambda, \sigma)$, where $2m = 
\alpha (C_j)$. Let
$n_0,n_1, \dotsc, n_l$ be a sequence of integers 
satisfying  
$n_0 \ge n_1 \ge \cdots \ge n_l = 1 $, $n_0= 
\lambda$,$n_1=\sigma$,
$n_{i-1}+n_{i+1} \equiv 0 \; (\roman{mod} \; n_i)$, and 
$\sum_{i=0}^{l-1} 1/n_in_{i+1} = (1/2)\vert s(C_j) \vert $.
\par 
Let $\roman{Ch}(A_j)$ be the chorizo space shown in  Figure
3, which consists of a 
disk and $l+2$ spheres. Again we consider the spaces 
$\roman{Ch}(f^iA_j)\
(i=0,1, \dotsc,
 (m/2)-1) $ identical.

Now the generalized quotient $S_f$ is defined to be the 
union of $\roman{Ch}
(B)$ and
$\roman{Ch}(A_j)$\<'s , $A_j$ running over all the annuli 
in $A$. A natural 
projection $\pi :
 {\Sigma}_g \to S_g $ can be defined. \par

Let $C_{\pi}$ be the mapping cylinder of $\pi$. We 
construct an ``open
book'' $\overline M$ with a ``page'' $C_{\pi}$. (See [Ta, 
W].)
 Then $M= \roman{int} (\overline M)$ has a complex 
structure, and we obtain a
 degenerating family $(M,D, \varphi)$ whose monodromy 
coincides with $f$.
Blow down $(-1)$-curves, if any, in $M$. All of the 
process is topologically 
canonical, and we get the inverse map
$\sigma : {\Cal P}_g^{-} \to {\Cal S}_g $ of $\rho : {\Cal 
S}_g \to
 {\Cal P}_g^{-} $, proving its bijectivity.\par

\heading 4. Conjugacy invariants \endheading

We define the {\sl{partition graph}} $X_f$ associated with 
a pseudoperiodic
homeomorphism of negative twist $f:{\Sigma}_g \to 
{\Sigma}_g$ as follows: Let 
$C$ be a precise system. The vertices (resp. the edges) of 
$X_f$ are in 
one-to-one 
correspondence to the connected components $b$ of 
${\Sigma}_g-C$ (resp.
the curves $\{C_i\}$ in $C$). An edge $e(C_i)$ joins 
vertices $v(b)$ and
$v(b')$ if and only if 
$C_i$ is in the adherence of $b$ and also of $b'$. The {\sl{
refined partition graph}} ${\overline{X}}_f$ is obtained 
from $X_f$ by 
subdividing those edges $e(C_i)$ that correspond to 
amphidrome curves by
their middle points. \par

A periodic map ${\psi}_f:{\overline{X}}_f \to 
{\overline{X}}_f $ is 
induced from $f$. The quotient graph $Y_f = 
{\overline{X}}_f/{\psi}_f$ is
a weighted graph in the sense that each vertex (and each 
edge) carries a 
positive integer called the {\sl{weight}}, which is the 
number of the vertices
(resp. the edges) of ${\overline{X}}_f$ over the vertex 
(resp. the edge)
of $Y_f$.

\midspace{5.5pc}
\caption{{\smc Figure 3}}

The conjugacy class of the periodic action ${\psi}_f : 
{\overline{X}}_f \to
{\overline{X}}_f$ can be interpreted as a cohomology class 
$c_f$ in a 
suitably defined {\sl{weighted cohomology group}} 
$H_W^1(Y_f)$. The weighted
graph $Y_f$ also serves as the {\sl{decomposition 
diagram}} of $S_f$,
and there is a natural collapsing map ${\eta}_f :S_f \to 
Y_f$.

\proclaim{Theorem 2} The triple $(S_f,Y_f,c_f)$ determines 
the conjugacy class of the mapping class $[f\,]$. 
\endproclaim

Nielsen's set of invariants introduced in [Ni2] has 
exactly the same amount
of information as $(S_f,Y_f)$ but lacks $c_f$. Thus his 
assertion [Ni2, \S 15;
G, Theorem 13.4] that his invariants are complete 
is incorrect.
\par
Here is an example. Let ${\Sigma}_6$ be the surface shown 
in Figure
4, which
has a system of curves $C=\{ C_1,\dotsc,C_5\}$. Let $f_k : 
{\Sigma}_6 \to 
{\Sigma}_6\ (k=1,2)$ be a homeomorphism such that 
$f_k^5\vert ({\Sigma}_6-C)
\simeq \roman{id}$, $f_k(b_i)=b_{i+k}$,$f_k(C_j)=C_{j+k}$, 
and $s(C_j)=-1$.
(Indices are taken modulo 5.) Although Nielsen's 
invariants are the same for
$[f_1]$ and $[f_2]$, these mapping classes are not 
conjugate, 
because the actions on the partition graphs are not 
conjugate.

\topspace{16pc}
\caption{{\smc Figure 4}}

\heading 5. Acknowledgment  \endheading

The authors are grateful to Allan Edmonds, Takayuki Oda, 
and Hiroshige Shiga
 for their useful information and comments.

\enddocument